\documentclass[cmfonts,a4paper]{aipproc}
\layoutstyle{8d}
\usepackage{amssymb,amsmath,amsthm,verbatim}
\usepackage[mathcal]{euler}
 \font\palba=pplb at 20pt
 \font\palo=pplr at 14pt
 \font\paloa=pplr at 12pt
 \font\palom=pplr at 10 pt
\title{\palba Chernoff's bound forms}
\author{\palo M. Grendar, Jr. and \underline{M. Grendar}}
{
 address =
 {Institute of Measurement Science, Slovak Academy of Sciences (SAS) 
           D\'ubravsk\'a cesta 9, Bratislava, 842 19, Slovakia 
           \& Institute of Mathematics and Informatics of Mathematical Institute SAS
           and of Matej Bel University, 
           Bansk\'a Bystrica, Slovakia.
           umergren@savba.sk}
}

\def\vc{\mathbf}

\newtheorem{theorem}{Theorem}

\newtheorem{note}{Note}

\begin{document}
\begin{abstract}
Chernoff's bound  binds a tail probability (ie. $Pr(X \ge a)$, where $a \ge
EX$). Assuming that the distribution of $X$ is $Q$, the logarithm  of the
bound is known to be equal to the value of relative entropy (or minus
Kullback-Leibler distance) for $I$-projection $\hat P$ of $Q$ on a set
$\mathcal{H} \triangleq \{P: E_PX = a\}$. Here, Chernoff's bound is related
to Maximum Likelihood on exponential form and consequently implications for
the notion of complementarity  are discussed. Moreover, a novel form of the
bound is proposed, which expresses the value of the Chernoff's bound
directly in terms of the $I$-projection (or generalized $I$-projection).%
\end{abstract}

\maketitle

\section{\paloa Introduction}

Originally developed as an asymptotic result for partial sums of random
variables, Chernoff's bound \cite{Chernoff} was later recognized to be
valid 'for any $n$'. It permitted to formulate Chernoff's bound in the
following form
\begin{theorem}
 Let $X$ be a random variable such that 
 $E e^{\theta v(X)} < \infty$, for all $\theta\in\mathbf{R}$, where $v(X)$
 is a concave, non-decreasing function of $X$. Let $a\geq EX$,
 $a\in\mathbf{R}$.
 Then
\begin{subequations}
 \begin{align}
 \log P(X\geq a) &\leq \min_{\theta \in\mathbf{R}}\:
                 \log Ee^{\theta v(X)} - \theta v(a),\\
\intertext{or, equivalently}
 P(X\geq a) &\leq \min_{\theta \in\mathbf{R}} \: \frac{Ee^{\theta v(X)}}{e^{\theta v(a)}}
 \end{align}
\end{subequations}

\end{theorem}

Since a proof of the Theorem (see for instance \cite{Weiss}) will be used
in the sequel, it will be recalled here.

\begin{proof}
Since $e^{\theta X}$ is a nonnegative valued and monotone function of $X$,
for $\theta > 0$ it is increasing in $X$. By assumption $v(X)$ is a
 non-decreasing function of X. Thus, by Markov's inequality
\begin{multline*}
P(X\ge a) = P(\theta X \ge \theta a) = P(\theta v(X) \ge \theta v(a)) = \\
= P\left(e^{\theta v(X)} \ge e^{\theta v(a)}\right) \le \frac{E e^{\theta
v(X)}}{e^{\theta v(a)}}
\end{multline*}
The inequality holds trivially for $\theta = 0$, thus the tightest bound is
achieved by minimizing the right-hand side expression, over $\theta \ge 0$.

To show that
\begin{equation*}{}
\arg\, \min_{\theta \ge 0}\frac{E e^{\theta v(X)}}{e^{\theta v(a)}} \equiv
\arg\, \min_{\theta \in \mathbf{R}}\frac{E e^{\theta v(X)}}{e^{\theta v(a)}}
\end{equation*}
apply Jensen's inequality both  to the exponential function and to
$v(\cdot)$, then recall that $a \ge E X$ and consequently realize, that
point of minimum of $E e^{\theta(v(X)
- v(a))}$ should occur for non-negative value of $\theta$.%

Hence,
\begin{equation*}
P(X\ge a) \le \min_{\theta \in \mathbf R} \frac{E e^{\theta
v(X)}}{e^{\theta v(a)}}
\end{equation*}

\end{proof}

Notation: Let us denote
\begin{equation}\label{eq:teta}
\hat\theta \triangleq \arg\: \min_{\theta \in\mathbf{R}} \: \log Ee^{\theta
v(X)} - \theta v(a)
\end{equation}
The entire right-hand side of (1a), (1b) will be denoted $C(a, v(\cdot),
\hat\theta)$, $C^P(a, v(\cdot), \hat\theta)$, respectively.

While it may appear at  first glance surprising, Chernoff's bound on tail
probability for a single random variable can be expressed in terms of
quantities related to a random sample of asymptotic size.  This is recalled
and summarized in the next two sections. The last, relatively self-standing
section, introduces a novel form/interpretation of Chernoff's bound.

\section{\paloa Chernoff's bound as a minimum  of I-divergence}

In this and the next section it will be  assumed that $X$ is either a
continuous random variable with pdf $g(X)$ defined on a support $S$; or a
discrete random variable with an $m$-element pmf $\vc q$.

First, the continuous case. Let $\mathcal{H}$ denote a class of pdf's,
$\mathcal{H} \triangleq \{f: E_f v(X) = v(a)\}$. Consider the following {\it
I-divergence minimization task} which consists of selecting  a pdf $\hat
f(X)$ from the class $\mathcal{H}$ that is closest to $g(X)$, where the
closeness is measured by $I$-divergence (or {\it I-distance})
$$
 I(f\parallel g) = E_f \log\frac{f(X)}{g(X)}
$$

Employing calculus of variations, it is possible to  show (see for instance
\cite{minimax}) that the unique solution (in open form) of the above task is
$$
\hat f(x) = \frac{g(x) e^{\hat\theta v(x)}}{\int_S g(x) e^{\hat\theta
v(x)}} = g(x) e^{\hat\theta v(x) - \log E_g e^{\hat\theta v(X)}}
$$
where $\hat\theta$ is a solution of
\begin{equation}\label{eq:constr}
E_{\hat f} v(X) = v(a)
\end{equation}

Consequently, it can be easily seen that the value of the $I$-divergence for
the pdf $\hat f(x)$ closest to $g(x)$ at the class $\mathcal{H}$ is
$$
I(\hat f\parallel g\mid f \in \mathcal{H}) = \hat\theta v(a) - \log E_g
e^{\hat\theta v(X)}
$$

Recalling the convex analysis duality theorem (see for instance
\cite{Ellis}), it can be shown that $\hat\theta$ which solves
(\ref{eq:constr}) and $\hat\theta$ of (\ref{eq:teta}) are the same.

Thus,
\begin{equation}
C(a, v(\cdot), \hat\theta) = - I(\hat f\parallel g\mid f \in \mathcal{H})
\end{equation}

In words, the logarithm of  tail probability of obtaining a value greater
than $a$ is bounded by the negative of the value of the $I$-distance of pdf
$\hat f(x)$ closest to $g(x)$ in the class $\mathcal{H}$ of all pdf's with
value of $E_{f} v(X)$ just equal to  $v(a)$.

Equivalent to the $I$-divergence minimization task is a {\it
relative-entropy maximization task} (since relative entropy $H(f\parallel
g) \triangleq - I(f\parallel g)$), thus
\begin{equation}\label{eq:rem}
C(a, v(\cdot), \hat\theta) = H(\hat f\parallel g\mid f \in \mathcal{H})
\end{equation}

The discrete case allows for deeper reading. Let now $\mathcal{H}$ denote a
class of pmf's, $\mathcal{H} \triangleq \{\vc p: E_p v(X) = v(a)\}$. The
relative entropy maximization (REM) task
\begin{equation}
\arg\, \max_{\vc p \in \mathcal{H}}\,\, - \sum_{i=1}^m p_i \log
\left(\frac{p_i}{q_i}\right)
\end{equation}
is solved by $\hat p_i = q_i e^{\hat\theta v(x_i) - \log E_q e^{\hat\theta
v(X)}}$, where $\hat\theta$ solves $E_{\hat p} v(X) = v(a)$. Consequently,
arguing along the same line as in the continuous case leads to the
conclusion similar to (\ref{eq:rem}),
\begin{equation}
C(a, v(\cdot), \hat\theta) = H(\hat{\vc p}\parallel{\vc q}\mid \vc{p} \in
\mathcal{H})
\end{equation}
which can now be  followed further to get
\begin{equation}\label{eq:cbrem}
C^P(a, v(\cdot), \hat\theta) = \prod_{i=1}^m {\left(\frac{q_i}{\hat
p_i}\right)}^{\hat p_i}
\end{equation}

Recalling the MaxProb justification of REM (see \cite{GG}) it can be noted
that $\hat{\vc p}$ is a limit of sequence of the most probable occurrence
vectors; and this way Chernoff's bound becomes related to random sample of
asymptotic size.

{\bf Example}. Let $X$ be defined on support $[1\ 2\ 3\ \dots\ 8]$ with pmf
$\vc q = [0.05\ 0.4\ 0.2\ 0.15\ 0.10$ $0.07\ 0.02\ 0.01]$. Thus $EX = 3.19$.
Setting $a = 4$ we ask for tail probability $P(X\ge 4)$ which is obviously
$0.35$. The closest in $I$-divergence to $\vc q$ pmf can be found to be
$\hat{\vc p} = [0.0236\ 0.2526\ 0.1692\ 0.1699\ 0.1517\ 0.1422$ $0.0544\
0.0364]$. Chernoff's bound calculated by (\ref{eq:cbrem}) then gives the
value $0.8829$. For $a=5$ it gives 0.5675, as compared to true $0.2$; for
$a=6$ it gives $0.27$, (true value is  0.1); and for $a=7$ it gives $0.087$
(true value is $0.03$).

\section{\paloa Chernoff's bound and Maximum Likelihood}

Let us assume a random sample $\vc X = \vc x$  of size $n$, such that
\begin{equation*}
\frac{1}{n}\sum_{i=1}^n v(x_i) = v(a)
\end{equation*}
where $a$, $v(\cdot)$ are given.

Let the supposed population from which the sample came be of the following
exponential form
$$
p_i(\theta) = q_i e^{\theta v(x_i) - \log E_q e^{\theta v(X)}}
$$
where $\vc q$ is a pmf, thus $\vc p$ is the exponentially tilted $\vc q$.

Maximum likelihood (ML) task lays in searching out a value of $\theta$
which is the most likely to generate the sample $\vc x$. The ML estimator
$\theta_{ML}$ of $\theta$ is known to be the solution of the likelihood
equation which is now just
$$
v(a) = \frac{E_q v(X) e^{\theta v(X)}}{E_q  e^{\theta v(X)}}
$$
Thus, $\hat\theta_{ML}\equiv \hat\theta$ (see also \cite{minimax}).

It is then interesting to relate Chernoff's bound to the above ML task. The
log-likelihood
$$
l(\theta) = \sum_{i=1}^m n_i\left(\log q_i + \theta v(x_i) -
\log\sum_{i=1}^m q_i e^{\theta v(x_i)}\right)
$$
where $n_i$ is occurrence of the $i$-th element of support at the sample.
So,
\begin{subequations}
\begin{align}{}
\frac{l(\hat\theta_{ML})}{n} &= \sum_{i=1}^m \frac{n_i}{n} \log q_i - C(a,
v(\cdot), \hat\theta), \\
\intertext{or equivalently, with $L_{\hat\theta}$ denoting the likelihood at
maximum,}
C^P(\cdot) &= \sqrt[n]{\dfrac{\prod_{i=1}^m q_i^{n_i}}{L_{\hat\theta}}}
\label{eq:mlc}
\end{align}
\end{subequations}
which establish ML-Chernoff's bound links.

Do they? For instance (\ref{eq:mlc}), combined with  (\ref{eq:cbrem}), lead
to conclusion
\begin{equation*}
\prod_{i=1}^m {\left(\frac{q_i}{\hat p_i}\right)}^{\hat p_i} =
\prod_{i=1}^m {\left(\frac{q_i}{\hat p_i}\right)}^{\frac{n_i}{n}}
\end{equation*}
which is false, except for the case when $\hat p_i \equiv \frac{n_i}{n}, i =
1, 2, \dots, m$\footnote{And except for the trivial case $q_i/\hat p_i =
1/m$, for all $i$}. This case happens to appear just for the random sample
of asymptotic size. Which solves the contradiction: since REM is indeed the
method which operates with a random sample of infinite size (c.f. \cite{GG},
or \cite{Vasicek}).

ML and REM tasks are complementary, regardless of sample size (see
\cite{minimax}). But, as the above 'deduction' shows, objective functions of
 both tasks (maximum likelihood, relative entropy, respectively) attain
a compatible relationship only when infinite sample size is assumed. And
this is indeed the case, because REM requires assumption about infiniteness
of random sample.

At the asymptotic, thanks to a conditional weak law of large numbers (see
 \cite{Vasicek}), Chernoff's bound is linked to the exponential
 form Maximum Likelihood by
$$
\frac{l(\hat\theta_{ML})}{n} \xrightarrow{p} \sum_{i=1}^m \hat p_i \log q_i
- C(\cdot)
$$
which leads further to the conclusion (similar in spirit to the Asymptotic
E\-qui\-par\-ti\-tion Property)
$$
\frac{l(\hat\theta_{ML})}{n} \xrightarrow{p} - H(\hat{\vc p} |\vc p \in
\mathcal H)
$$
where $H(\vc p) \triangleq - \sum p_i \log p_i$ is Shannon's entropy.

\section{\paloa New form of Chernoff's bound}

The logarithm of the tail probability $\log Pr(X \ge a)$ cannot exceed the
convex conjugate of the cumulant generating function, of the random variable
$v(X)$ --- this is a statement of the 'log-Chernoff bound' (recall (1a)),
for the log-tail-probability. Assuming that the distribution of $X$ is $Q$,
the value of the log-Chernoff's bound becomes equal to negative of the
value of the Kullback-Leibler distance ($I$-divergence) for $I$-projection
$\hat P$ of $Q$ on a set $\mathcal H \triangleq \{P: E_P v(X) = v(a)\}$,
recall (4). Under the assumption, the Chernoff's bound value can also be
expressed directly in terms of $I$-projection -- as  will be shown here.

In order to make it relatively self-standing and precise, it will be given
in terms of measure theory and $I$-projection (see \cite{Csiszar}). Though
the presented variant of Chernoff's bound is the same in the case of a
discrete random variable as well as  in the case of a continuous one, each
case will be discussed under different existence considerations, hence its
formulation is separated into separate theorems.

\subsection{\paloa Discrete measure}

\begin{theorem}
 Let $(\Omega, \mathcal{F}, Q)$ be a countable probability space and let
 $X: \Omega \rightarrow \mathbf{R}$ be a random variable taking values
 $\{x_1, x_2, \dots\}$. Let $a \in \mathbf{R}$ such that $a \ge E_Q X$.
 Assume that $E_Q e^{\theta X} < \infty$ for all $\theta \in \mathbf R$.
 Let $\mathcal{P}$ denote the class of all probability measures on
 $(\Omega, \mathcal{F})$ and $\mathcal{H} = \{P \in \mathcal{P}: E_P X = a\}$.
 If $a$ is in the convex hull of $\{x_1, x_2, \dots\}$ $\mathcal{H} \neq\emptyset$.
 Assume this to be the case. Let $\hat P$ be the $I$-projection of $Q$ on
 $\mathcal{H}$, that is $I(\hat{P}\| Q) = \inf_{P \in \mathcal{H}} I(P\| Q)$.
 If $I(\hat P\| Q)$ is finite, then
 $$
  Q(\omega: X(\omega) \ge a) \le \frac{Q(a)}{\hat P(a)}
 $$ %
\end{theorem}

\begin{proof}
  To save space, let $p_i \triangleq P\{x_i\}$, $q_i \triangleq Q\{x_i\}$,
  $\hat p_i \triangleq \hat P\{x_i\}$.

  By (\cite{Ellis}, Thm II.5.2, Thm VIII.3.1), under the assumptions,
  the $I$-projection of $Q$ on $\mathcal H$ exists, it is unique, and has
  the following form
 \begin{equation}\label{eq:expo}
 \hat p_i = q_i e^{\hat\lambda x_i - \log E_Q e^{\hat\lambda X}}
 \end{equation}
 where
 \begin{equation}\label{eq:equiv}
 \hat\lambda = \arg\, \min_{\lambda} E_Q e^{\lambda(X - a)}
 \end{equation}
 exists and it is unique.

 Since $a \ge EX$, $E_Q e^{\theta X} < \infty$
 for all $\theta \in \mathbf R$, the standard
 proof of Chernoff's bound (see the Introduction) guarantees that
 $$
  \min_{\theta \in R} E_Q e^{\theta(X-a)} \ge Q(X \ge a)
 $$
 or, with use of (\ref{eq:equiv})
 \begin{equation}\label{eq:CB}
  E_Q e^{\hat\lambda(X-a)} \ge Q(X \ge a)
 \end{equation}

  Noting that $\hat{P}(a) = Q(a) e^{\hat\lambda a - \log E_Q e^{\hat\lambda X}}$
  then shows that the LHS of (\ref{eq:CB}) is just $\frac{Q(a)}{\hat{P}(a)}$,
  which completes the proof.
\end{proof}

\begin{note}
 The claim of Theorem 2 could be directly extended by replacing $X$ by
 any concave, non-decreasing and bounded function $v(X)$.
\end{note}

\subsection{\paloa Absolutely continuous measure}

 Let now a measurable function $X: \Omega \rightarrow \mathbf R$, defined on
 a probability space $(\Omega, \mathcal{F}, \mu)$ induces on $\mathbf R$
 a law $Q$ dominated by Lebesgue measure $\lambda$,
 so that its density $q(X)$ with respect to $\lambda$ exists.
 Let $\mathcal H$ be a convex set of
 laws $P$ on $\mathbf R$ whose densities $p(X)$ with respect to Lebesgue measure
 exist.
 $I$-projection  $\hat P$ of $Q$ on $\mathcal H$ is then such
 $\hat P \in \mathcal H$ that
 $I(\hat P\| Q) = \inf_{P \in \mathcal H} I(P\|Q)$, where
 $I(P\| Q) = \int p(x) \log \frac{p(x)}{q(x)} \lambda(dx)$. There,
 $0 \log 0 = 0$, $\log \frac{b}{0} = + \infty$
 conventions are assumed\footnote{The definition
 of $I$-projection was adapted from  \cite{Csiszar}. Throughout the
 paper $\log$ denotes the natural logarithm (though it is in fact immaterial
 for the claims which are made).}.

 Assuming existence of $I$-projection, the new form of Chernoff's bound
 can be stated as follows:

 \begin{theorem}
  Let $v(X)$ be a concave and non-decreasing function of $X$.
  Let $a \ge E_Q X$,
  $a \in \mathbf R$. Let $\mathcal H \triangleq \{p: E_P v(X) = v(a)\}$.
  Let $\hat p(x)$ -- the density
  corresponding to the I-projection of $Q$ on $\mathcal H$ -- exist.
  Let $E_Q e^{\theta v(X)} < \infty$, $E_Q v(X) e^{\theta v(X)} < \infty$,
  for all $\theta \in \mathbf R$. Then
  $$
   \mu(\omega: X(\omega) \ge a) \le \frac{q(a)}{\hat p(a)}
  $$
  provided that $q(a) \ne 0$, $\hat p(a) \ne 0$ and that $a$ is the point where both
  $\hat{p}(X)$ and $q(X)$ are unique.
 \end{theorem}

 \begin{proof}
  By Theorem 3.1 and Corollary 3.1 of \cite{Csiszar}
  $I$-projection of $Q$ on $\mathcal H$ has a density with respect
  to Lebesgue measure of the following open form
  $\hat p(x, \eta) = e^{\eta v(x) - \log E_Q e^{\eta v(X)}}$,
  which is closed by $\hat\eta$ such that $E_{\hat P} v(X) = v(a)$.
  The density is unique, up to a set $\aleph$ of measure zero.

  By assumptions $E_Q e^{\theta v(X)} \le \infty, \forall \theta \in \mathbf R$,
  $E_Q v(X) e^{\theta v(X)} \le \infty, \forall \theta \in \mathbf R$
  so $\hat\theta \triangleq \arg \min_{\theta \in \mathbf R} E_Q e^{\theta (v(X) - v(a))}$
  exists and it is unique. The assumptions also guarantee that differentiation of
  $E_Q e^{\theta (v(X) - v(a))}$ with respect to $\theta$ can be performed
  under integral (cf. \cite{Durrett}, Theorem A(9.1)). Consequently, it can
  be directly seen that $\hat\theta$ solves $E_{\hat P} v(X) = v(a)$ and is
  identical with $\hat\eta$.

  (The above argument could be also made by invoking
  (\cite{Ellis}, Thm VIII.3.1).)

  It is assumed that $\hat p \notin \aleph$, and different than zero as is also
  assumed $q(a)$, thus
  \begin{equation}\label{eq:link}
  E_Q e^{\hat\theta(v(X) - v(a))} = \frac{q(a)}{\hat p(a)}
  \end{equation}

  The assumption $E_Q e^{\theta v(x)} < \infty$, $\forall\theta\in\mathbf R$
  together with assumed properties of $v(\cdot)$
  guarantee validity of Chernoff's bound claim:
  \begin{equation}\label{eq:CCB}
    E_Q e^{\hat\theta(v(X) - v(a))} \ge \mu(X\ge a)
  \end{equation}

  Comparing (\ref{eq:link}) and (\ref{eq:CCB}) completes the proof.

 \end{proof}

As far as the existence of $I$-projection is concerned, Csisz\'ar's work
(see \cite{Csiszar}, discussion on pp. 151, 154 and Theorems 2.1, 3.2)
implies that for the case considered above, if $I(P\| Q) < \infty$ for some
$P \in \mathcal H$ and if $\mathcal H \ne \emptyset$ and if $v(X)$ is
bounded then the $I$-projection $\hat P$  of $Q$ on $\mathcal H$ exists, it
is unique, and has the form  $\hat p(x) = q(x) e^{\hat\theta v(x) - \log
E_Q e^{\hat\theta v(X)}}$.

Though the $I$-projection may not exist in the case of unbounded $v(X)$,
nevertheless {\it generalized $I$-projection} introduced by Tops\o e (see
\cite{Flemming}) and studied further by Csisz\'ar (see \cite{Csiszar2})
exists and take up the exponential form, which -- even in this case --
permits to formulate Chernoff's bound  in terms of generalized
$I$-projection. This will be done after a brief reminder of generalized
$I$-projection, which is adapted from \cite{Csiszar2}.

 Let $(S, \mathcal{B})$ be a measurable space, $X$ -- random variable, and $P$,$Q$ be two
probability measures defined on the measurable space. $I$-divergence
$I(P||Q)$ between them is

\begin{displaymath}
 I(P||Q) = \left\{ \begin{array}{ll}
 \int \log(dP/dQ) \, dP & \textrm{if $P \ll Q$} \\
 +\infty & \textrm{otherwise}
 \end{array} \right.
\end{displaymath}
and let $\mathcal H$ be a set of probability measures on $(S,\mathcal B)$.
Let
$$
I(\mathcal H||Q) \triangleq \inf_{P\in\mathcal H} I(P||Q)
$$

 Generalized $I$-projection $\hat P$ of $Q$ on $\mathcal H$ is such a
probability measure not necessarily in $\mathcal H$ that every sequence of
probability measures $P_n \in \mathcal H$ with $I(P_n||Q) \rightarrow
I(\mathcal H||Q)$ converges to $\hat P$ in variation.

Making use of Csisz\'ar's results, the generalized $I$-projection form of
Chernoff's bound can be stated as follows:

 \begin{theorem}
  Let $v(X)$ be a concave, non-decreasing, not necessarily bounded function of $X$.
  Let $a \ge E_Q X$,
  $a \in \mathbf R$.
  Let $E_Q e^{\theta (v(X) - v(a))}$ attain its minimum at $\hat\theta$.
  Let $\mathcal H \triangleq \{P: E_P v(X) = v(a)\}$.
  Let $\frac{d\hat P}{dQ}(x)$ be the generalized $I$-projection of $Q$
  on $\mathcal H$.

  Then
  $$
     Pr(X\ge a) \le  \frac{1}{\frac{d\hat P}{dQ} (a)}
  $$
 provided that unique $\frac{d\hat P}{dQ} (a) \neq 0$.
 \end{theorem}

 \begin{proof}
 Since $\hat\theta$ exists (by assumption), by  (cf. \cite{Csiszar2}, p. 778)
 the generalized $I$-projection of $Q$ on $\mathcal H$ is
 \begin{equation}\label{eq:GIP}
  \frac{d\hat P}{dQ} (x) = e^{\hat\theta v(X) - \log E_Q e^{\hat\theta v(X)}}
 \end{equation}
 Thus, $1/\frac{d\hat P}{dQ} (a)$
 is just $E_Q e^{\hat\theta(v(X) - v(a))}$,
 ie. the Chernoff's bound value, which binds $Pr(X\ge a)$.
 \end{proof}

\section{\paloa Acknowledgments}

Hospitality of Banach Centre of the Institute of Mathematics, Polish
Academy of Sciences, where a part of this study was performed as a part  of
the European Community Center of Excellence programme (package 'Information
Theory and its Applications to Physics, Finance and Biology') is gratefully
acknowledged. The work was also supported by the grant VEGA 1/7295/20 from
the Scientific Grant Agency of the Slovak Republic.
 It is a pleasure to thank an anonymous reviewer, Chris Williams,
 Franti\v sek Rubl\'\i k and Viktor Witkovsk\'y
 for  valuable comments which helped to  improve both the contents and form of
 this work.

\renewcommand\refname

\medskip\medskip

M. Grendar, Jr. and M. Grendar, ``Chernoff's bound forms,'' in {\it
Bayesian inference and Maximum Entropy methods in Science and Engineering,}
edited by Ch. Williams, AIP Conference Proceedings 659, Melville, New York,
2003, pp. 67-72.

\bigskip\medskip
\rightline{\palom to Mar}
\rightline{\palom June 23, 2003}
\end{document}